\documentclass{ifacconf}

\usepackage{graphicx}      
\usepackage{natbib}        
\usepackage{amsmath} 
\usepackage{amssymb}  
\usepackage{xcolor}
\usepackage{mathtools}
\usepackage{enumerate}
\usepackage{url}

\newtheorem{Theorem}{Theorem}
\newtheorem{Lemma}{Lemma}
\newtheorem{Problem}{Problem}
\newtheorem{Remark}{Remark}
\newtheorem{Corollary}{Corollary}

\newtheorem{Assumption}{Assumption}
\newtheorem{Definition}{Definition}
\newcommand{\bb}{\boldsymbol}
\DeclareMathOperator{\R}{\mathbb R}
\DeclareMathOperator*{\argmin}{arg\,min}

\begin{document}
\begin{frontmatter}

\title{Safe Control Design for Unknown Nonlinear Systems with Koopman-based Fixed-Time Identification
} 


\author[First]{Mitchell Black} 
\author[Second]{Dimitra Panagou} 

\address[First]{Department of Aerospace Engineering,  University of Michigan, Ann Arbor, MI 48109, USA (e-mail: mblackjr@umich.edu).}
\address[Second]{Department of Robotics and Department of Aerospace Engineering,  University of Michigan, Ann Arbor, MI 48109, USA (e-mail: dpanagou@umich.edu)}

\begin{abstract}
We consider the problem of safe control design for a class of nonlinear, control-affine systems subject to an unknown, additive, nonlinear disturbance. Leveraging recent advancements in the application of Koopman operator theory to the field of system identification and control, we introduce a novel fixed-time identification scheme for the infinitesimal generator of the infinite-dimensional, but notably linear, Koopman dynamical system analogous to the nonlinear system of interest. That is, we derive a parameter adaptation law that allows us to recover the unknown, residual nonlinear dynamics in the system within a finite-time independent of an initial estimate. We then use properties of fixed-time stability to derive an error bound on the residual vector field estimation error as an explicit function of time, which allows us to synthesize a provably safe controller using control barrier function based methods. We conduct a quadrotor-inspired case study in support of our proposed method, in which we show that safe trajectory tracking is achieved despite unknown, nonlinear dynamics.
\end{abstract}


\end{frontmatter}


\section{Introduction}

Recent advances in computing power and memory storage have ushered in an era of estimation, identification, and control for autonomous systems dominated by data-driven methods. For example, compared to the 74 kilobytes of memory available on the United States National Aeronautics and Space Administration's (NASA) first lunar module computer, the gigabytes of memory used in many of today's data-driven approaches to dynamical system identification (e.g. deep neural networks) have allowed engineers to create significantly more expressive models. Though regression methods are widely-used for linear system identification, the field of identification for nonlinear systems is vast. Popular approaches in recent years include classes of neural networks (NNs), including deep NNs (e.g. \cite{zancato2021novel}) and recurrent NNs for time-varying systems (\cite{gonzalez2018non}), Gaussian processes (\cite{frigola2013integrated}), and more recently the application of Koopman operator theory (e.g. \cite{Mauroy2020Koopman,brunton2016koopman,Klus2020KoopmanGenerator}, among others), which introduces an infinite-dimensional but notably linear representation of a nonlinear system on which traditional linear identification approaches may be used.

Under Koopman theory there exists a linear Koopman dynamical system that captures the dynamics of the original nonlinear system over an infinite-dimensional space of scalar functions known as observables. Beginning with \cite{Mauroy2020Koopman}, recent work has focused on using data-driven approaches to approximate a finite-dimensional matrix representation of the Koopman operator, which acts as a state-transition operator for the Koopman dynamical system. In particular, extended dynamic mode decomposition (EDMD), first introduced in \cite{williams2015data} has emerged as a popular tool for carrying out such an approximation. The end result in many cases is a batch estimate of either the Koopman matrix (i.e. in \cite{Bruder2021DataDriven,Haseli2021Approximation}) or its infinitesimal generator (\cite{Klus2020KoopmanGenerator,Drmac2021KoopmanGenerator}) obtained by solving a least-squares regression problem. Potential shortcomings of this class of approaches include slower response times than e.g. recursive methods, and a lack of formal guarantees on the approximation error bound, which may be particularly detrimental when used in control design. In contrast, it has been shown by \cite{Black2022Fixed} that fixed-time stability in the context of recursive parameter identification admits a such bound on the identification error as an explicit function of time.

Finite- and fixed-time stability (FTS and FxTS) are stronger notions of stability for equilibria of a dynamical system, each of which guarantee convergence of the system trajectories to the origin within a finite time. They have been used in the analysis of linear parameter identification schemes by \cite{rios2017time,Ortega2022Parameter}, and synthesized for the purpose of safe control design in \cite{Black2022Fixed,Wang2022Robust}. The benefit to recursive parameter identification in fixed-time, i.e. in a finite-time independent of the initial condition, is the knowledge of an error bound on the identification error as an explicit function of time. When synthesized with a safe control law, this class of identification schemes yields less conservative control solutions, as highlighted in \cite{Black2022Fixed}.

Control barrier functions (CBFs) have proven to be a useful tool for safe control synthesis. As a model-based approach, however, it is critical that an accurate system model be available in order to preserve forward invariance of the set of safe states. Though robust CBF controllers can protect against bounded disturbances to the system dynamics (e.g. \cite{Jankovic2018Robust, Black2020Quadratic}), the cost is conservatism. Various other approaches to safe control have sought to adapt to the unknown residual dynamics (e.g. \cite{Taylor2020Adaptive, Lopez2020Robust}), or to learn their effects via data-driven Koopman-based policies both online (\cite{folkestad2020data}) and offline (\cite{zinage2022neural}). None of these methods, however, provide guarantees on learning convergence time. 

In this paper, we address this open problem by introducing a Koopman-based identification scheme for safe control design that guarantees convergence within a fixed-time for a class of nonlinear, control-affine systems subject to an additive, nonlinear perturbation. We use knowledge of the bound on convergence time to quantify the identification error as an explicit function of time, the magnitude of which is leveraged to design a provably safe CBF-based controller. We demonstrate the advantages of our proposed approach on a trajectory tracking problem, and highlight that the identification and control laws succeed in preserving safety of the system even in the presence of measurement noise.






The rest of the paper is organized as follows. In Section \ref{sec.preliminaries} we introduce the preliminaries and define the problem under consideration. Section \ref{sec.fxt_estimation} contains our main result on fixed-time nonlinear system identification, which we use in Section \ref{sec.robust_adaptive_control} to design a safe controller. We demonstrate the approach on a numerical case study in Section \ref{sec.case_studies}, and conclude in Section \ref{sec.conclusion} with directions for future work.

\section{Preliminaries and Problem Statement}\label{sec.preliminaries}

In this paper, we use the following notation. $\mathbb R$ denotes the set of real numbers. The ones matrix of size $n \times m$ is denoted $\bb{1}_{n \times m}$. We use $\|\cdot\|$ to denote the Euclidean norm and $\|\cdot\|_{\infty}$ to denote the sup norm. We denote the minimum and maximum eigenvalue of a matrix $\boldsymbol{M}$ as $\lambda_{min}(\boldsymbol{M})$ and $\lambda_{max}(\boldsymbol{M})$, and refer to its r$^{th}$ singular value as $\sigma_r(\bb{M})$, to its nullspace as $\mathcal{N}(\boldsymbol{M})$, and its i$^{th}$ column as $\mathrm{col}_i(\bb{M})$. The gradient operator is $\nabla$, and the Lie derivative of a function $V:\mathbb R^n\rightarrow \mathbb R$ along a vector field $f:\mathbb R^n\rightarrow\mathbb R^n$ at a point $x\in \mathbb R^n$ is denoted as $L_fV(x) \triangleq \frac{\partial V}{\partial x} f(x)$.

Consider the following class of nonlinear, control-affine systems
\begin{equation}\label{eq.nonlinear_system}
    \dot{\bb{x}} = f(\bb{x}(t)) + g(\bb{x}(t))\bb{u}(t) + d(\bb{x}(t)), \quad \bb{x}(0) = \bb{x}_0,
\end{equation}
where $\bb{x} \in \mathcal{X} \subset \R^n$ and $\bb{u} \in \R^m$ denote the state and control input vectors, the drift vector field $f: \R^n \rightarrow \R^n$ and control matrix field $g: \R^n \rightarrow \R^n \times \R^m$ are known and continuous, and $d: \R^n \rightarrow \R^n$ is an unknown disturbance known to be continuous and to obey $\|d(\bb{x})\|_{\infty} \leq D < \infty$ for all $\bb{x} \in \mathcal{X}$.
Consider also the following set of safe states,
\begin{equation}\label{eq.safe_set}
    S = \{\bb{x} \in \mathcal{X} \; | \; h(\bb{x}) \geq 0\},
\end{equation}
for a continuously differentiable function $h: \R^n \rightarrow \R$, where the boundary and interior of $S$ are $\partial S = \{\bb{x} \in \R^n \; | \; h(\bb{x}) = 0\}$ and $\textrm{int}(S) = \{\bb{x} \in \R^n \; | \; h(\bb{x}) > 0\}$ respectively. The trajectories of \eqref{eq.nonlinear_system} are said to be \textit{safe} if the set $S$ is \textit{forward-invariant}, i.e. if $\bb{x}_0 \in S \implies \bb{x}(t) \in S, \forall t \geq 0$. The following lemma, known as Nagumo's Theorem, provides necessary and sufficient conditions for rendering $S$ forward-invariant.
\begin{Lemma}(\cite{blanchini1999set})\label{lem.nagumos_theorem}
    Suppose that $\bb{u}(t)$ is continuous such that the closed-loop trajectories of \eqref{eq.nonlinear_system} are uniquely determined in forward-time. The set $S$ is forward-invariant if and only if
    \begin{equation}\label{eq.forward_invariance_condition}
        \dot{h} = \frac{\partial h(\bb{x})}{\partial \bb{x}}\dot{\bb{x}
} \geq 0, \; \forall \bb{x} \in \partial S.
    \end{equation}
\end{Lemma}
In recent years, control barrier functions have emerged as a viable approach for control design satisfying \eqref{eq.forward_invariance_condition}.
\begin{Definition}(\cite{ames2017control})\label{def.cbf}
    Given a set $S \subseteq \mathcal{X} \subset \R^n$ defined by \eqref{eq.safe_set} for a continuously differentiable function $h: \R^n \rightarrow \R$, the function $h$ is a \textbf{control barrier function} (CBF) defined on the set $\mathcal{X}$ if there exists a Lipschitz continuous class $\mathcal{K}_\infty$ function $\alpha: \R \rightarrow \R$ such that
    \begin{equation}\label{eq.cbf_condition}
        \sup_{\bb{u} \in \R^m}\dot{h}(\bb{x}, \bb{u}) \geq -\alpha(h(\bb{x})),
    \end{equation}
    for all $\bb{x} \in \mathcal{X}$.
\end{Definition}
We refer to \eqref{eq.cbf_condition} as the CBF condition, and observe that it constitutes sufficiency for the satisfaction of \eqref{eq.forward_invariance_condition}. As such, any continuous control law $\bb{u}(t)$ that 1) admits unique closed-loop trajectories of \eqref{eq.nonlinear_system} in forward-time and 2) satisfies \eqref{eq.cbf_condition} renders the trajectories of \eqref{eq.nonlinear_system} safe. Consider now that for the system \eqref{eq.nonlinear_system} the CBF condition is
\begin{equation}
    \sup_{\bb{u} \in \R^m}\left[L_fh(\bb{x}) + L_gh(\bb{x})\bb{u} + L_dh(\bb{x})\right] \geq -\alpha(h(\bb{x})), \nonumber
\end{equation}
where, without identification of $d(\bb{x})$, the precise value of $L_dh(\bb{x})$ is unknown. It is known, however, that
\begin{equation}
    -b_d \leq L_{d}h(\bb{x}) \leq b_d, \nonumber
\end{equation}
where $b_d = D \left|\frac{\partial h(\bb{x})}{\partial \bb{x}}\right| \mathbf{1}_{n \times 1}$. Under such circumstances, a robust-CBF may be used for safe control design.
\begin{Definition}(\cite{Jankovic2018Robust})\label{eq.robust_cbf}
    Given a set $S \subseteq \mathcal{X} \subset \R^n$ defined by \eqref{eq.safe_set} for a continuously differentiable function $h: \R^n \rightarrow \R$, the function $h$ is a \textbf{robust control barrier function} (r-CBF) for the system \eqref{eq.nonlinear_system} defined on the set $\mathcal{X}$ if there exists a Lipschitz continuous class $\mathcal{K}_\infty$ function $\alpha: \R \rightarrow \R$ such that
    \begin{equation}\label{eq.rcbf_condition}
        \sup_{\bb{u} \in \R^m}\left[L_fh(\bb{x}) + L_gh(\bb{x})\bb{u} -b_d\right] \geq -\alpha(h(\bb{x})),
    \end{equation}
    for all $\bb{x} \in \mathcal{X}$.
\end{Definition}

Designing a controller to protect against the worst possible disturbance in perpetuity, however, may lead to poor performance, especially if $D$ is large. Recent work (e.g. \cite{Lopez2020Robust,Black2022Fixed}) has shown that this may be mitigated by using an estimate of the unknown disturbance $\hat d(\bb{x})$. Thus, we define the vector field estimation error $\Tilde{d}(\bb{x})$ as
\begin{equation}\nonumber
    \Tilde{d}(\bb{x}) \coloneqq d(\bb{x}) - \hat d(\bb{x}).
\end{equation}
In \cite{Black2022Fixed}, it was shown under mild assumptions that if the uncertain vector field is parameter-affine, i.e. if 
\begin{equation}\nonumber
    d(\bb{x}) = \Delta(\bb{x})\bb{\theta}^*,
\end{equation}
for some known, continuous, bounded regressor matrix $\Delta: \mathcal{X} \rightarrow \R^{n \times p}$ and unknown, static, polytopic parameters $\bb{\theta}^* \in \Theta \subset \R^p$, then the vector field estimation error may be driven to zero within a fixed time using parameter adaptation, i.e. $\|\Delta(\bb{x}(t))(\bb{\theta}^* - \hat{\bb{\theta}}(t))\| \rightarrow 0$ as $t \rightarrow T < \infty$, independent of $\hat{\bb{\theta}}(0)$. We now review the notion of fixed-time stability.

\subsection{Fixed-Time Parameter Identification}
Consider a nonlinear, autonomous system of the form
\begin{equation}\label{eq.autonomous_nonlinear_system}
    \dot{\bb{x}} = F(\bb{x}), \quad \bb{x}(0) = \bb{x}_0,
\end{equation}
where $F: \R^n \rightarrow \R^n$ is continuous such that \eqref{eq.autonomous_nonlinear_system} admits a unique solution for all $\bb{x}_0 \in \R^n$, the value of which at time $t$ is denoted $\bb{\varphi}_t(\bb{x}_0)$, and where $F(0) = 0$.
\begin{Definition}(\cite{polyakov2012nonlinear})\label{def.fxts}
    The origin of \eqref{eq.autonomous_nonlinear_system} is fixed-time stable (FxTS) if it is stable in the sense of Lyapunov and any solution $\bb{\varphi}_t(\bb{x}_0)$ of \eqref{eq.autonomous_nonlinear_system} reaches the origin within a finite time $T$ independent of $\bb{x}_0$, i.e. $\exists T < \infty$ such that $\bb{\varphi}_t(\bb{x}_0) = 0$ for all $t \geq T$, $\forall \bb{x}_0 \in \R^n$.
\end{Definition}
In what follows, we review a fixed-time stable parameter adaptation law from the literature.
\begin{Theorem}(\cite{Black2022Fixed})\label{thm.fixed_time_parameter_adaptation}
    Consider a perturbed dynamical system of the form \eqref{eq.nonlinear_system}. Suppose that the following hold:
    \begin{enumerate}[i)]   
        \item the unknown, additive dynamics are parameter-affine, i.e. $d(\bb{x}) = \Delta(\bb{x})\bb{\theta}^*$,
        \item there exist a known matrix $\bb{M}(t) \in \R^{n \times p}$ and vector $\bb{v}(t) \in \R^n$ such that $\bb{M}(t)(\bb{\theta}^* - \hat{\bb{\theta}}(t)) = \bb{v}(t)$,
        \item the nullspace of $\Delta(\bb{x}(t))$ is constant for all $t \leq T$, i.e. $\mathcal{N}(\Delta(\bb{x}(t))) = \mathcal{N}(\Delta(\bb{x}(0)))$, $\forall t \leq T$,
    \end{enumerate}
    where
    \begin{equation}\label{eq: FxT Time Bound}
        T = \frac{\mu\pi}{2k_V^2\sqrt{ab}},
    \end{equation}
    with $a,b>0$, $\mu > 2$, and
    \begin{equation}\label{eq: kV}
        k_V = \sigma_r(\boldsymbol{M})\sqrt{2\lambda_{max}(\boldsymbol{\Gamma})},
    \end{equation}
    where $\boldsymbol\Gamma \in \mathbb R^{p\times p}$ is a constant, positive-definite, gain matrix and $\sigma_r(\boldsymbol{M})>0$ denotes the smallest nonzero singular value of $\boldsymbol{M}$ over the time interval.
    Then, under the ensuing parameter adaptation law,
    \begin{equation}\label{eq.adaptation_law}
        \dot{\hat{\boldsymbol{\theta}}} = \boldsymbol\Gamma \boldsymbol{M}^T\boldsymbol{v}\left(a\|\boldsymbol{v}\|^{ \frac{2}{\mu}} + \frac{b}{\|\boldsymbol{v}\|^{\frac{2}{\mu}}}\right),
    \end{equation}
    the estimated disturbance $\hat d(\bb{x}(t))$ converges to the true disturbance $d(\bb{x}(t))$ within fixed-time $T$, i.e. $\Delta(\bb{x}(t))\hat{\bb{\theta}}(t) \rightarrow \Delta(\bb{x}(t))\bb{\theta}^*$ as $t \rightarrow T$, and $\Delta(\bb{x}(t))\hat{\bb{\theta}}(t) = \Delta(\bb{x}(t))\bb{\theta}^*$ for all $t \geq T$, independent of $\hat{\bb{\theta}}(0)$.
\end{Theorem}
\begin{pf}
    See \cite[Proof of Theorem 3]{Black2022Fixed}.
\end{pf}

Theorem \ref{thm.fixed_time_parameter_adaptation} provides a framework for adapting parameter estimates $\hat{\bb{\theta}}$ such that an unknown disturbance of the form $d(\bb{x}) = \Delta(\bb{x})\bb{\theta}^*$ is learned within fixed-time. In reality, however, it is far more common for the unknown vector field $d(\bb{x})$ to be nonlinear, which to this point has precluded the use of \eqref{eq.adaptation_law} as a learning or adaptation strategy. By utilizing Koopman operator theory, however, we can transform the problem of identifying the nonlinear function $d$ into a linear, albeit infinite-dimensional, identification problem, which with appropriate modifications permits the use of the above adaptation framework. 

\subsection{Koopman Operator based Identification}
Koopman theory dictates that a nonlinear system of the form \eqref{eq.autonomous_nonlinear_system} has an analogous and notably linear representation in an infinite-dimensional Hilbert space $\mathcal{Q}$ consisting of continuous, real-valued functions $q: \mathcal{X} \rightarrow \R$ referred to as \textit{observables}. The continuous-time Koopman dynamical system analogous to \eqref{eq.autonomous_nonlinear_system} is then described by
\begin{equation}\label{eq.lifted_dynamics}
    \dot q = \mathcal{L}q, \quad q \in \mathcal{Q},
\end{equation}
where $\mathcal{L}$ denotes the infinitesimal generator of the linear semigroup of Koopman operators $\mathcal{U}^t: \mathcal{Q} \rightarrow \mathcal{Q}$, i.e.
\begin{equation}\nonumber
    \mathcal{L}q = \lim_{t \rightarrow 0}\frac{\mathcal{U}^tq - q}{t} = F \cdot \nabla q.
\end{equation}

For tractability, however, many works (e.g. \cite{Bruder2021DataDriven,Drmac2021KoopmanGenerator}, among others) derive matrix representations $\bb{U} \in \R^{N \times N}$ and $\bb{L} \in \R^{N \times N}$ of the respective finite-rank operators $\mathcal{U}_N^t = \Pi_N\mathcal{U}^t|_{\mathcal{Q}_N}$ and $\mathcal{L}_N = \Pi_N\mathcal{L}|_{\mathcal{Q}_N}$, where $\Pi_N: \mathcal{Q} \rightarrow \mathcal{Q}_N$ is a projection operator onto the subspace $\mathcal{Q}_N \subset \mathcal{Q}$ (spanned by $N>n$ linearly independent basis functions $\{\psi_i: \mathcal{X} \rightarrow \R\}_{i=1}^N$) and $\mathcal{O}|_{\mathcal{Q}_N}$ denotes the restriction of the operator $\mathcal{O}$ to $\mathcal{Q}_N$. We refer the reader to \cite{mauroy2020koopmanbook} for additional details, and instead highlight that in practice $\bb{U}$ and $\bb{L}$ are taken to be the respective solutions to
\begin{align}
    \bb{\psi}^T(\bb{x})\bb{U} = (\bb{\psi}(\bb{\varphi}_t(\bb{x})))^T, \label{eq.koopman_matrix_regression} \\
    \bb{L}^T \bb{\psi}(\bb{x}) = \frac{\partial \bb{\psi}(\bb{x})}{\partial \bb{x}} F(\bb{x}), \label{eq.koopman_generator_regression}
\end{align}
where $\bb{\psi}(\bb{x}) = [\psi_1(\bb{x}) \hdots  \psi_N(\bb{x})]^T \in \R^N$ and $\frac{\partial \bb{\psi}(\bb{x})}{\partial \bb{x}} \in \R^{N \times n}$.

If $\bb{L}$ can be identified directly (as in e.g. \cite{Klus2020KoopmanGenerator}), the vector field $F$ may be reconstructed by solving \eqref{eq.koopman_generator_regression} for $F(\bb{x})$. When this is not possible, identification of $\bb{U}$ may be used to reconstruct $F$ after computing $\bb{L}$ via
\begin{equation}\label{eq.koopman_generator_log}
    \bb{L} = \frac{1}{T_s}\log \bb{U},
\end{equation}
in the case of sampled data, where $\log$ denotes the principal matrix logarithm and $T_s>0$ is the sampling interval. We observe that both \eqref{eq.koopman_matrix_regression} and \eqref{eq.koopman_generator_regression} describe linear systems of equations of the form $\bb{a}^T\bb{X} = \bb{b}$, and thus $\bb{X}$ (in this case $\bb{U}$ or $\bb{L}$) can be identified using linear identification techniques such as the parameter identification law \eqref{eq.adaptation_law}.

\subsection{Problem Statement}
Now, reconsider the unknown, control-affine, nonlinear system \eqref{eq.nonlinear_system}. Suppose that an estimate of its Koopman generator matrix $\hat{\bb{L}}$ is available, and let the estimated unknown vector field $\hat d(\bb{x})$ then via \eqref{eq.koopman_generator_regression} be the solution to
\begin{equation}\nonumber
    \hat{\bb{L}}^T\bb{\psi}(\bb{x}) = \frac{\partial \bb{\psi}(\bb{x})}{\partial \bb{x}} \big(f(\bb{x}) + g(\bb{x})\bb{u} + \hat d(\bb{x})\big).
\end{equation}
We assume that $\frac{\partial \bb{\psi}(\bb{x})}{\partial \bb{x}}$ is full column rank, which may be satisfied by design (e.g. sinusoidal basis functions), and thus have that $\hat d(\bb{x}) \rightarrow d(\bb{x})$ as $\hat{\bb{L}} \rightarrow \bb{L}$ (which can also be satisfied if $\hat{\bb{U}} \rightarrow \bb{U}$). Define the vectorized Koopman matrix and generator ($\bb{\mu}^*$ and $\bb{\lambda}^*$), and their estimates ($\hat{\bb{\mu}}$ and $\hat{\bb{\lambda}}$), as 
\begin{align}
    \bb{\mu}^* &\coloneqq [\mathrm{col}_1^T(\bb{U}) \hdots \mathrm{col}_N^T(\bb{U})]^T, \label{eq.koopman_matrix_vector} \\
    \bb{\lambda}^* &\coloneqq [\mathrm{col}_1^T(\bb{L}) \hdots \mathrm{col}_N^T(\bb{L})]^T, \label{eq.koopman_generator_vector} \\
    \hat{\bb{\mu}} &\coloneqq [\mathrm{col}_1^T(\hat{\bb{U}}) \hdots \mathrm{col}_N^T(\hat{\bb{U}})]^T, \label{eq.koopman_matrix_estimate} \\
    \hat{\bb{\lambda}} &\coloneqq [\mathrm{col}_1^T(\hat{\bb{L}}) \hdots \mathrm{col}_N^T(\hat{\bb{L}})]^T,\label{eq.koopman_generator_estimate}
\end{align}
and observe that for the system \eqref{eq.nonlinear_system} the relations \eqref{eq.koopman_matrix_regression} and \eqref{eq.koopman_generator_regression} are equivalent to
\begin{equation}\label{eq.matrix_vector_regression}
    \bb{\Psi}(\bb{x})\bb{\mu}^* = (\bb{\psi}(\bb{\varphi}_t(\bb{x})))^T,
\end{equation}
and
\begin{equation}\label{eq.generator_vector_regression}
    \bb{\Psi}(\bb{x})\bb{\lambda}^* = \frac{\partial \bb{\psi}(\bb{x})}{\partial \bb{x}}\big(f(\bb{x}) + g(\bb{x})\bb{u} + d(\bb{x})\big),
\end{equation}
respectively, where
\begin{equation}\label{eq.psi_matrix}
    \bb{\Psi}(\bb{x}) \coloneqq \begin{bmatrix}
        \bb{\psi}^T(\bb{x}) & 0 & \hdots & 0 \\
        0 & \bb{\psi}^T(\bb{x}) & \hdots & 0 \\
        \vdots & & \ddots & \vdots \\ 
        0 & \hdots & 0 & \bb{\psi}^T(\bb{x})
    \end{bmatrix} \in \R^{N \times N^2}.
\end{equation}

Let the Koopman matrix and Koopman generator estimation errors respectively be denoted
\begin{align}
    \Tilde{\bb{\mu}} &= \bb{\mu}^* - \hat{\bb{\mu}}, \nonumber \\
    \Tilde{\bb{\lambda}} &= \bb{\lambda}^* - \hat{\bb{\lambda}}, \nonumber
\end{align}
and observe that $\bb{\Psi}(\bb{x})\hat{\bb{\lambda}} = \bb{\Psi}(\bb{x})\bb{\lambda}^*$ for all $\Tilde{\bb{\lambda}} \in \mathcal{N}(\bb{\Psi}(\bb{x}))$. 

We are now ready to formally define the problem under consideration.
\begin{Problem}\label{prob.main_problem}
    Consider a dynamical system of the form \eqref{eq.nonlinear_system}. Design adaptation and control laws, $\dot{\hat{\bb{\lambda}}} = \eta(\bb{x}, \bb{u}, \hat{\bb{\lambda}})$ and $\bb{u} = \kappa(\bb{x}, \hat{\bb{\lambda}})$ respectively, such that
    \begin{enumerate}
        \item the Koopman generator error vector, $\Tilde{\bb{\lambda}}$, is rendered fixed-time stable to the nullspace of $\bb{\Psi}(\bb{x})$, i.e. $\Tilde{\bb{\lambda}}(t) \rightarrow \mathcal{N}(\bb{\Psi}(\bb{x}))$ as $t \rightarrow T$ and $\Tilde{\bb{\lambda}}(t) \in \mathcal{N}(\bb{\Psi}(\bb{x}))$ for all $t \geq T$, independent of $\hat{\bb{\lambda}}(0)$, and
        \item the system trajectories remain safe for all time, i.e. $\bb{x}(t) \in S$,  $\forall t \geq 0$.
    \end{enumerate}
\end{Problem}
In the ensuing section, we introduce our approach to solving the first element of Problem \ref{prob.main_problem}.


\section{Nonlinear Estimation in Fixed-Time}\label{sec.fxt_estimation}
In this section, we introduce our proposed adaptation law $\dot{\hat{\bb{\lambda}}} = \eta(\bb{x}, \bb{u}, \hat{\bb{\lambda}})$ for the fixed-time identification of the Koopman generator vector $\bb{\lambda}$, which allows us to identify the unknown vector field $d(\bb{x})$ in \eqref{eq.nonlinear_system} within a fixed-time. Before introducing one of our main results, we require the following assumptions.
\begin{Assumption}\label{ass.koopman_projection}
    The projection of the infinite-dimensional Koopman operator $\mathcal{U}^t$ onto the finite-rank subspace $\mathcal{Q}_N$ exactly describes the evolution of observables $q \in \mathcal{Q}$, i.e. $\mathcal{U}_N^t q = (\Pi_N\mathcal{U}^t)q$, for all $q \in \mathcal{Q}$. 
\end{Assumption}
\begin{Assumption}\label{ass.rank_Psi}
    There exist scalars $s>0$, $T>0$ such that $\sigma_N(\bb{\Psi}(\bb{x}(t))) \geq s$ for all $0 \leq t \leq T$, where $\bb{\Psi}(\bb{x}(t))$ is given by \eqref{eq.psi_matrix}.
\end{Assumption}
The satisfaction of Assumption \ref{ass.koopman_projection} depends on the choice of $N$ (and thus on the basis functions $\bb{\psi}$), and while generally this is an open problem recent work has studied the existence of Koopman invariant subspaces (see e.g. \cite{brunton2016koopman}), i.e. subspaces $\mathcal{Q}_N \subset \mathcal{Q}$ over which Assumption \ref{ass.koopman_projection} holds. For our numerical study in Section \ref{sec.case_studies}, we find that bases $\bb{\psi}$ constructed using 
monomials or sinusoids
work well.
The satisfaction of Assumption \ref{ass.rank_Psi} evidently depends on the choice of basis functions $\psi_i$. Note, however, that $\bb{\Psi}(\bb{x}(t))$ is guaranteed to be full row-rank (which implies that $\sigma_N(\bb{\Psi}(\bb{x}(t))) > 0$) provided that $\exists i \in [N]$ such that $\psi_i(\bb{x}(t)) \neq 0$. This can be guaranteed with an appropriate choice of bases, e.g. $\psi_1(\bb{x}(t)) = 1$.

\begin{Theorem}\label{thm.koopman_adaptation}
    Suppose that Assumptions \ref{ass.koopman_projection} and \ref{ass.rank_Psi} hold, where
    \begin{equation}\label{eq.koopman_generator_settling_time}
        T = \frac{w\pi}{4s\lambda_{max}(\bb{\Gamma})\sqrt{ab}},
    \end{equation}
    with $a, b>0$, $w > 2$, and $\bb{\Gamma} \in \R^{N^2 \times N^2}$ a constant, positive-definite gain matrix. Then, under the ensuing adaptation law
    \begin{equation}\label{eq.koopman_adaptation_law}
    \begin{aligned}\small
        \dot{\hat{\bb{\lambda}}} = \bb{\Gamma}\bb{\Psi}^T(\bb{x})\bb{\nu}(\bb{x},\hat{\bb{\lambda}})\left(a\|\bb{\nu}(\bb{x},\hat{\bb{\lambda}})\|^{2 / w} + \frac{b}{\|\bb{\nu}(\bb{x},\hat{\bb{\lambda}})\|^{2 / w}}\right),
    \end{aligned}
    \end{equation}
    the Koopman generator error vector $\Tilde{\bb{\lambda}}$ is rendered FxTS to the nullspace of $\bb{\Psi}(\bb{x})$, i.e. $\Tilde{\bb{\lambda}}(t) \rightarrow \mathcal{N}(\bb{\Psi}(\bb{x}(t)))$ as $t \rightarrow T$ and $\Tilde{\bb{\lambda}}(t) \in \mathcal{N}(\bb{\Psi}(\bb{x}))$ for all $t \geq T$, independent of $\hat{\bb{\lambda}}(0)$, where
    \begin{equation}\label{eq.nu_adaptation_vector}
        \bb{\nu}(\bb{x},\hat{\bb{\lambda}}) = \frac{\partial \bb{\psi}(\bb{x})}{\partial \bb{x}}\dot{\bb{x}} - \bb{\Psi}(\bb{x})\hat{\bb{\lambda}}.
    \end{equation}
\end{Theorem}
\begin{pf}
    We first show that there exists a time-invariant Koopman generator vector $\bb{\lambda}(t) = \bb{\lambda}^*$, $\forall t \geq 0$, and then prove that under \eqref{eq.koopman_adaptation_law} the aassociated Koopman generator error vector $\Tilde{\bb{\lambda}}$ is rendered FxTS to $\mathcal{N}(\bb{\Psi}(\bb{x}))$.
    
    First, under Assumption \ref{ass.koopman_projection} it follows that there exists a finite-rank operator $\mathcal{L}_N: \mathcal{Q}_N \rightarrow \mathcal{Q}_N$ such that the nonlinear dynamics of \eqref{eq.nonlinear_system} may be represented by the following linear system in the space of observables:
    \begin{equation}\nonumber
        \dot q = \mathcal{L}_Nq, \quad q \in \mathcal{Q}.
    \end{equation}
    Then, there exists a finite-dimensional matrix representation $\bb{L} \in \R^{N \times N}$ in a basis $\{\psi_i: \mathcal{X} \rightarrow \R\}_{i=1}^N$ corresponding to the operator $\mathcal{L}_N$ such that the relation \eqref{eq.koopman_generator_regression} holds over the trajectories of \eqref{eq.nonlinear_system}. Thus, the Koopman generator matrix $\bb{L}$ admits the (time-invariant) Koopman generator vector $\bb{\lambda}^*$ defined by \eqref{eq.koopman_generator_vector}.

    Next, observe that \eqref{eq.generator_vector_regression} over the trajectories of \eqref{eq.nonlinear_system} may be modified to obtain
    \begin{align}
        \bb{\Psi}(\bb{x})\bb{\lambda}^* - \bb{\Psi}(\bb{x})\hat{\bb{\lambda}} &= \frac{\partial \bb{\psi}(\bb{x})}{\partial \bb{x}}\dot{\bb{x}} - \bb{\Psi}(\bb{x})\hat{\bb{\lambda}}, \nonumber \\
        \bb{\Psi}(\bb{x})\Tilde{\bb{\lambda}} &= \bb{\nu}(\bb{x},\hat{\bb{\lambda}}), \nonumber
    \end{align}
    where $\bb{\nu}(\bb{x},\hat{\bb{\lambda}})$ is given by \eqref{eq.nu_adaptation_vector}. Thus, we have that the premises of Theorem \ref{thm.fixed_time_parameter_adaptation} are satisfied with $\bb{M} = \bb{\Psi}$ and $\bb{v} = \bb{\nu}$ and the adaptation law \eqref{eq.koopman_adaptation_law} takes the form of \eqref{eq.adaptation_law}. Then, with Assumption \ref{ass.rank_Psi} it follows directly from Theorem \ref{thm.fixed_time_parameter_adaptation} that $\Tilde{\bb{\lambda}}$ is rendered FxTS to $\mathcal{N}(\bb{\Psi}(\bb{x}))$ with settling time given by \eqref{eq.koopman_generator_settling_time}. 
\end{pf}

In what follows, we show how the parameter adaptation law \eqref{eq.koopman_adaptation_law} results in learning the exact disturbance $d(\bb{x})$ to the system dynamics \eqref{eq.nonlinear_system} within fixed-time.
\begin{Corollary}\label{cor.disturbance_fixed_time}
    Consider the system \eqref{eq.nonlinear_system}. Suppose that the premises of Theorem \ref{thm.koopman_adaptation} hold, and that the estimated Koopman vector $\hat{\bb{\lambda}}$ is adapted according to \eqref{eq.koopman_adaptation_law}. If the estimated disturbance $\hat d(\bb{x})$ is taken to be
    \begin{equation}\label{eq.dhat}
        \hat d(\bb{x}(t)) = \frac{\partial \bb{\psi}(\bb{x}(t))}{\partial \bb{x}}^{\dagger}\bb{\Psi}(\bb{x}(t))\hat{\bb{\lambda}}(t) - a(\bb{x}(t),\bb{u}(t)),
    \end{equation}
    where $a(\bb{x}(t), \bb{u}(t)) = f(\bb{x}(t)) + g(\bb{x}(t))\bb{u}(t)$, then,
    the vector field estimation error $\Tilde{d}(\bb{x}(t))$ is rendered FxTS to the origin and the estimated disturbance $\hat{d}(\bb{x}(t))$ converges to the true disturbance $d(\bb{x}(t))$ within a fixed-time $T$ given by \eqref{eq.koopman_generator_settling_time}, i.e. $\Tilde{d}(\bb{x}(t)) \rightarrow 0$ and $\hat d(\bb{x}(t)) \rightarrow d(\bb{x}(t))$ as $t \rightarrow T$ independent of $\hat{d}(\bb{x}(0))$.
\end{Corollary}
\begin{pf}
    We first observe from \eqref{eq.generator_vector_regression} that the disturbance $d(\bb{x}(t))$ is the solution to
    \begin{equation}\label{eq.d_regression}
    \begin{aligned}\small
        \frac{\partial \bb{\psi}(\bb{x}(t))}{\partial \bb{x}}d(\bb{x}(t)) = \bb{\Psi}(\bb{x}(t))\bb{\lambda}^* - \frac{\partial \bb{\psi}(\bb{x}(t))}{\partial \bb{x}}a(\bb{x}(t),\bb{u}(t)).
    \end{aligned}\normalsize
    \end{equation}
    Next, it follows from Theorem \ref{thm.koopman_adaptation} that under \eqref{eq.koopman_adaptation_law} $\hat{\bb{\lambda}}(t) \rightarrow \bb{\lambda}^*$ as $t \rightarrow T$. Then, we have that $\bb{\Psi}(\bb{x}(t))\hat{\bb{\lambda}}(t) \rightarrow \bb{\Psi}(\bb{x}(t))\bb{\lambda}^*$ and thus that $\frac{\partial \bb{\psi}(\bb{x}(t))}{\partial \bb{x}}\hat d(\bb{x}(t)) \rightarrow \frac{\partial \bb{\psi}(\bb{x}(t))}{\partial \bb{x}}d(\bb{x}(t))$ as $t \rightarrow T$ when $\hat d(\bb{x}(t))$ is taken to be the solution to \eqref{eq.d_regression}.
    Finally, with $\frac{\partial \bb{\psi}(\bb{x}(t))}{\partial \bb{x}}$ full column rank we use its pseudoinverse $\frac{\partial \bb{\psi}(\bb{x}(t))}{\partial \bb{x}}^{\dagger}$ to recover \eqref{eq.dhat} and thus have that $\hat d(\bb{x}(t)) \rightarrow d(\bb{x}(t))$ as $t \rightarrow T$.

\end{pf}

For the purpose of control design it is important to know how the estimation error signals behave during the transient period $t \leq T$ before the unknown vector field $d(\bb{x})$ has been learned. In contrast to least-squares and related regression based approaches to learning the Koopman matrix $\bb{U}$ and/or generator matrix $\bb{L}$, our FxTS parameter adaptation law allows us to derive explicit estimation error bounds as a function of time. In fact, prior work (see \cite{Black2022Fixed}) has shown that the magnitude of this error bound is a monotonically decreasing function of time. In the following result, we introduce a modification to the prior work in order to derive a bound on the magnitude of the vector field estimation error $\Tilde{d}(\bb{x}(t))$ as an explicit function of time.
\begin{Corollary}\label{cor.error_bound}
    Suppose that the premises of Corollary \ref{cor.disturbance_fixed_time} hold. If, in addition, the initial estimated Koopman generator vector is set to zero, i.e. $\hat{\bb{\lambda}}(0) = \mathbf{0}_{N^2 \times 1}$, and $\boldsymbol{\Gamma}$ in \eqref{eq.koopman_adaptation_law} is constant, positive-definite, and also \textit{diagonal}, then $\forall t \in [0,T]$, where $T$ is given by \eqref{eq.koopman_generator_settling_time}, the following expression constitutes a monotonically decreasing upper bound on $\|\Tilde{d}(\bb{x}(t))\|_{\infty}$:
    \begin{equation}\label{eq.err_bound}
        \|\Tilde{d}(\bb{x}(t))\|_{\infty} \leq \Lambda \sigma_{max}(\bb{W}(t)) \tan^{\frac{w}{2}}(A(t)) \coloneqq \delta(t),
    \end{equation}
    where
    \begin{equation}\label{eq.Lambda_gain}
        \Lambda = \sqrt{2 \lambda_{max}(\boldsymbol{\Gamma})}\left(\frac{a}{b}\right)^{w/4},
    \end{equation}
    and
    \begin{align}
        \bb{W}(t) &= \frac{\partial \bb{\psi}(\bb{x}(t))}{\partial \bb{x}}^{\dagger}\bb{\Psi}(\bb{x}), \label{eq.w_mat} \\
        A(t)&=\max\left\{\Xi - \frac{\sqrt{ab}}{w}t,0\right\}, \label{eq.A(t)} \\ 
        \Xi &= \tan^{-1}\left(\sqrt{\frac{b}{a}}\left(\frac{1}{2}\bb{l}^T\Gamma^{-1}\bb{l}\right)^{\frac{1}{w}}\right), \label{eq.xi_term}
    \end{align}
    where $\bb{l} = \frac{2D}{\sigma_{min}(\bb{W}(0))} \cdot \boldsymbol{1}_{N^2\times1}$, and $\|\Tilde{d}(\bb{x}(t))\|_{\infty} =0$, $\forall t > T$.
\end{Corollary}
\begin{pf}
    See Appendix \ref{app.proof_d_error}.
\end{pf}

Knowledge of the upper bound on the disturbance estimation error bound \eqref{eq.err_bound} permits the use of robust, adaptive model-based control techniques. In particular, we will show in the next section how to synthesize a CBF-based controller that guarantees safety both before and after the transient phase $t \leq T$ during which the unknown disturbance $d(\bb{x})$ is learned, and in doing so address the second element of Problem \ref{prob.main_problem}.

\section{Robust-Adaptive Control Design}\label{sec.robust_adaptive_control}
In this section, we describe two approaches to synthesizing the Koopman-based parameter adaptation law with a CBF-based control law for safe control under model uncertainty.

\subsection{Robust-CBF Approach}
In the first approach, we demonstrate how to apply robust-CBF principles to the design of a safe controller $\bb{u} = \kappa(\bb{x},\hat{\bb{\lambda}})$ when using the Koopman-based adaptation scheme \eqref{eq.koopman_adaptation_law}.

\begin{Theorem}
    Consider a system of the form \eqref{eq.nonlinear_system}, a safe set $S$ defined by \eqref{eq.safe_set} for a continuously differentiable function $h: \mathcal{X} \rightarrow \R$, and suppose that the premises of Corollary \ref{cor.error_bound} hold. Then, any control input $\bb{u}$ satisfying
    \begin{equation}\label{eq.robust_koopman_cbf_condition}
        \sup_{\bb{u} \in \R^m}\left[L_fh(\bb{x}) + L_gh(\bb{x})\bb{u} + L_{\hat d}h(\bb{x}) - b_d(t)\right] \geq -\alpha(h(\bb{x}))
    \end{equation}
    renders the trajectories of \eqref{eq.nonlinear_system} safe, where
    \begin{equation}\label{eq.bd_bound}
        b_d(t) = \left|\frac{\partial h}{\partial \bb{x}}\right|\delta(t) \cdot \mathbf{1}_{n \times 1},
    \end{equation}
    and $\delta(t)$ is given by \eqref{eq.err_bound}.
\end{Theorem}
\begin{pf}
    Observe that over the trajectories of \eqref{eq.nonlinear_system}
    \begin{align}
        \dot h &= L_fh(\bb{x}) + L_gh(\bb{x})\bb{u} + L_dh(\bb{x}) \nonumber \\
        & = L_fh(\bb{x}) + L_gh(\bb{x})\bb{u} + \frac{\partial h}{\partial \bb{x}}\hat d(\bb{x}) + \frac{\partial h}{\partial \bb{x}}\Tilde d(\bb{x}) \nonumber \\
        &\geq L_fh(\bb{x}) + L_gh(\bb{x})\bb{u} + \frac{\partial h}{\partial \bb{x}}\hat d(\bb{x}) - \left|\frac{\partial h}{\partial \bb{x}}\right|\delta(t) \cdot \mathbf{1}_{n \times 1}. \nonumber
    \end{align}
    By Corollary \ref{cor.error_bound} it follows that $\|\Tilde d(\bb{x}(t))\|_{\infty} \leq \delta(t)$ for all $t \geq 0$. Therefore, $\dot h \geq -\alpha(h(\bb{x})$ whenever \eqref{eq.robust_koopman_cbf_condition} holds, and thus $S$ is rendered forward-invariant by any control input satisfying \eqref{eq.robust_koopman_cbf_condition}.
\end{pf}

It is worth noting that as the estimated disturbance $\hat d(\bb{x})$ converges to the true disturbance $d(\bb{x})$ the robustness term $b_d(t)$ will go to zero. So while initially the condition \eqref{eq.robust_koopman_cbf_condition} may demand large control inputs to guarantee safety in the face of a the unknown disturbance, as $t \rightarrow T$ the term $b_d(t) \rightarrow 0$ and the standard CBF condition is recovered.

\subsection{Robust-Adaptive CBF Approach}
In this approach, we define the following robust-adaptive safe set
\begin{equation}\label{eq.ra_safe_set}
    S_r = \{\bb{x} \in \mathcal{X}: h_r(\bb{x}, t) \geq 0\}
\end{equation}
for the continuously differentiable function
\begin{equation}\nonumber
    h_r(\bb{x}, t) = h(\bb{x}) - \frac{1}{2}\bb{\delta}^T(t)\bb{\Omega}^{-1}\bb{\delta}(t),
\end{equation}
for $\bb{\delta}(t) = \delta(t) \cdot \mathbf{1}_{n \times 1}$ with $\delta(t)$ given by \eqref{eq.err_bound}, and a constant, positive-definite matrix $\bb{\Omega} \in \R^{n \times n}$. We note that the set $S_r$ defined by \eqref{eq.ra_safe_set} is a subset of the safe set $S$ defined by \eqref{eq.safe_set}, i.e. $S_r \subseteq S$. We now introduce a robust-adaptive CBF condition that renders the trajectories of \eqref{eq.nonlinear_system} safe.

\begin{Theorem}
    Consider a system of the form \eqref{eq.nonlinear_system}, a set $S_r$ defined by \eqref{eq.ra_safe_set} for a continuously differentiable function $h_r: \mathcal{X} \rightarrow \R$, and suppose that the premises of Corollary \ref{cor.error_bound} hold. Then, any control input $\bb{u}$ satisfying
    \begin{equation}\label{eq.ra_koopman_cbf_condition}
        \sup_{\bb{u} \in \R^m}\left[L_fh_r(\bb{x}) + L_gh_r(\bb{x})\bb{u} - r\big(t, \hat d(\bb{x}(t))\big) \right] \geq -\alpha(h_r(\bb{x}))
    \end{equation}
    renders the trajectories of \eqref{eq.nonlinear_system} safe, where
    \begin{equation}\nonumber
        r\big(t, \hat d(\bb{x}(t))\big) = \mathrm{Tr}(\bb{\Omega}^{-1})\delta(t)\dot{\delta}(t) + b_d(t),
    \end{equation}
    where $\delta(t)$ is given by \eqref{eq.err_bound}, $b_d(t)$ is given by \eqref{eq.bd_bound}, and
    \begin{equation}\label{eq.delta_dot}
    \begin{aligned}
        \dot \delta(t) &= \Lambda\dot{\sigma}_{max}(\bb{W}(t))\tan^{\frac{w}{2}}(A(t)) \\
        &\quad - \frac{1}{2}\Lambda\sigma_{max}(\bb{W}(t))\sqrt{ab}\tan^{\frac{w}{2}-1}(A(t))\mathrm{sec}^2(A(t))
    \end{aligned}
    \end{equation}
\end{Theorem}
\begin{pf}
    Follows directly from \cite[Theorem 5]{Black2022Fixed} by replacing $\Tilde{\bb{\theta}}$ with $\Tilde{d}(\bb{x})$.
\end{pf}

\begin{Remark}
    We note that the robust-adaptive CBF condition \eqref{eq.ra_koopman_cbf_condition} requires the time-derivative of the maximum singular value of the matrix $\bb{W}(t)$ given by \eqref{eq.w_mat}, i.e. $\dot{\sigma}_{max}(\bb{W}(t))$. While this may not be available in closed-form, it may be approximated in practice using finite-difference methods.
\end{Remark}

Since both the robust \eqref{eq.robust_koopman_cbf_condition} and robust-adaptive \eqref{eq.ra_koopman_cbf_condition} CBF conditions ensure safety of the trajectories of \eqref{eq.nonlinear_system}, either condition may be included as an affine constraint in the now popular quadratic program based control law (eg. \cite{ames2017control,Black2020Quadratic}). We now introduce one such iteration of the QP controller,
\begin{subequations}\label{eq.cbf_qp_controller}
\begin{align}
    \bb{u}^* = \argmin_{\bb{u} \in \R^m} &\frac{1}{2}\|\bb{u}-\bb{u}^0\|^2 \label{subeq.cbf_qp_objective}\\
    \textrm{s.t.} \quad &\forall s\in[1..c] \nonumber \\
    \mathrm{Either} \; \eqref{eq.robust_koopman_cbf_condition} \; \mathrm{or} \; \eqref{eq.ra_koopman_cbf_condition}, \label{subeq.cbf_qp_constraints}
\end{align}
\end{subequations}
the objective \eqref{subeq.cbf_qp_objective} of which seeks to find a minimally deviating solution $\bb{u}^*$ from a nominal, potentially unsafe input $\bb{u}^0$ subject to the specified CBF constraint \eqref{subeq.cbf_qp_constraints}.

In the following section, we demonstrate the efficacy of our jointly proposed adaptation \eqref{eq.koopman_adaptation_law} and control \eqref{eq.cbf_qp_controller} laws on a quadrotor tracking problem.

\section{Numerical Case Study}\label{sec.case_studies}

Let $\mathcal{F}$ be an inertial frame with a point $s_0$ denoting its origin. Consider a quadrotor seeking to track a Gerono lemnisicate (i.e. figure-eight) trajectory amidst circular obstacles in the 2D plane. Quadrotor dynamics are known to be differentially-flat, thus as shown to be feasible in \cite{Zhou2014vector} we take the model to be the following 2D double-integrator subject to an unknown, wind disturbance:
\begin{equation}\label{eq.double_integrator}
    \begin{bmatrix}
        \dot x \\ \dot y \\ \dot{v}_x \\ \dot{v}_y
    \end{bmatrix} = \begin{bmatrix}
        v_x \\ v_y \\ a_x \\ a_y
    \end{bmatrix} + \begin{bmatrix}
        0 \\ 0 \\ d_x(\bb{z}) \\ d_y(\bb{z})
    \end{bmatrix},
\end{equation}
where $x$ and $y$ denote the position coordinates (in m), $v_x$ and $v_y$ are the velocities (in m/s), and $a_x$ and $a_y$ are the accelerations (in m/s$^2$). The full state and control input vectors are $\bb{z} = [x \; y \; v_x \; v_y]^T \in \R^4$ and $\bb{u} = [a_x \; a_y]^T \in \R^2$ respectively, and $d_x:\R^4 \rightarrow \R$ and $d_y:\R^4 \rightarrow \R$ are unknown wind-gust accelerations satisfying the requirements of $d$ in \eqref{eq.nonlinear_system}. Specifically, we used the wind-gust model from \cite{davoudi2020quad} to obtain spatially varying wind velocities $w_i(\bb{z})$ and set $d_i(\bb{z}) = C_d(w_i(\bb{z}) - v_i)$ for $i \in \{x, y\}$, where $C_d$ is a drag coefficient, such that $\|d_x(\bb{z})\|_{\infty}, \|d_y(\bb{z})\|_{\infty} \leq D = 10$.

We consider the presence of two circular obstacles, each of which occludes the desired quadrotor path. As such, the safe set is defined as
\begin{equation}\nonumber
    S = \{\bb{z} \in \R^4: h_1(\bb{z}) \geq 0\} \cap \{\bb{z} \in \R^4: h_2(\bb{z}) \geq 0\},
\end{equation}
where $h_i(\bb{z}) = (x - c_{x,i})^2 + (y - c_{y,i})^2 - R^2$ for $i \in \{1,2\}$, $(c_{x,i}, c_{y,i})$ denotes the center of the i$^{th}$ obstacle, and $R$ is its radius. Since $h_1,h_2$ are relative-degree two with respect to \eqref{eq.double_integrator}, we use future-focused CBFs for a form of safe, predictive control (see \cite{black2022intersection} for details).

We use forms of the CBF-QP control law\footnote{All simulation code and data are available online at \url{https://github.com/6lackmitchell/nonlinear-fxt-adaptation-control}} \eqref{eq.cbf_qp_controller} corresponding to both the robust \eqref{eq.robust_koopman_cbf_condition} and robust-adaptive \eqref{eq.ra_koopman_cbf_condition} CBF conditions, and compare the performance against a naive (i.e. assuming exact identification, $\hat d = d$) CBF controller equipped with the data-driven Koopman-based identification schemes proposed in \cite{Bruder2021DataDriven} and \cite{Klus2020KoopmanGenerator} respectively. For the robust and robust-adaptive simulations we inject additive Gaussian measurement noise into both $\bb{x}$ and $\dot{\bb{x}}$ in order to stress-test the algorithm under non-ideal conditions. We use the nominal control law introduced for quadrotors in \cite{Schoellig2012PeriodicQuadrotor} and adapted for our dynamics, where the reference trajectory is the Gerono lemniscate defined by
\begin{align}
    x^*(t) &= 4\sin(0.2\pi t)\nonumber \\
    y^*(t) &= 4\sin(0.2\pi t)\cos(0.2\pi t), \nonumber
\end{align}
which specifies that one figure-eight pattern be completed every 10s. 
Our circular obstacles are centered on $(-2.5, 0)$ and $(2, -1)$ respectively, each with a radius of $R=1.5$m. For all controllers, we used linear class $\mathcal{K}_\infty$ functions $\alpha(h) = h$. For our Koopman basis functions, we used sinusoids of the form $\psi_i = \sqrt{2}\cos(n\pi z)$, $\psi_{i+1} = \sqrt{2}\sin(n\pi z)$, for $n \in \{1,2\}$ and $z \in \{x, y, v_x, v_y\}$.

\begin{figure}
    \begin{center}
        \includegraphics[width=1\columnwidth,clip]{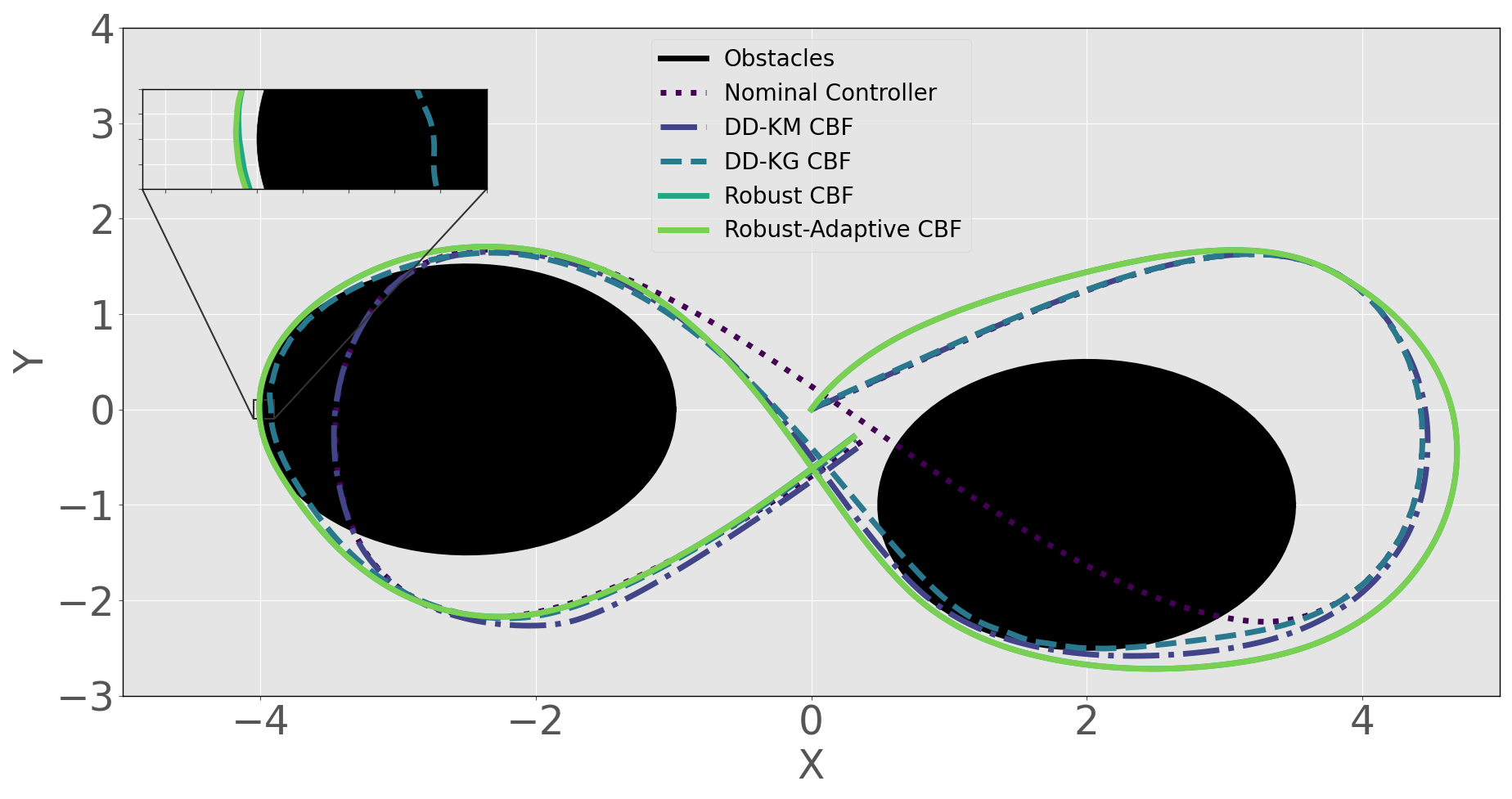}
        \caption{\small{XY paths under the various CBF-QP control laws in the double-integrator example. Only the controllers using the proposed Koopman-based fixed-time identification scheme succeed in preserving safety.}}
        \label{fig.xy_single_integrator}
    \end{center}
\end{figure}

\begin{figure}
    \begin{center}
        \includegraphics[width=1\columnwidth,clip]{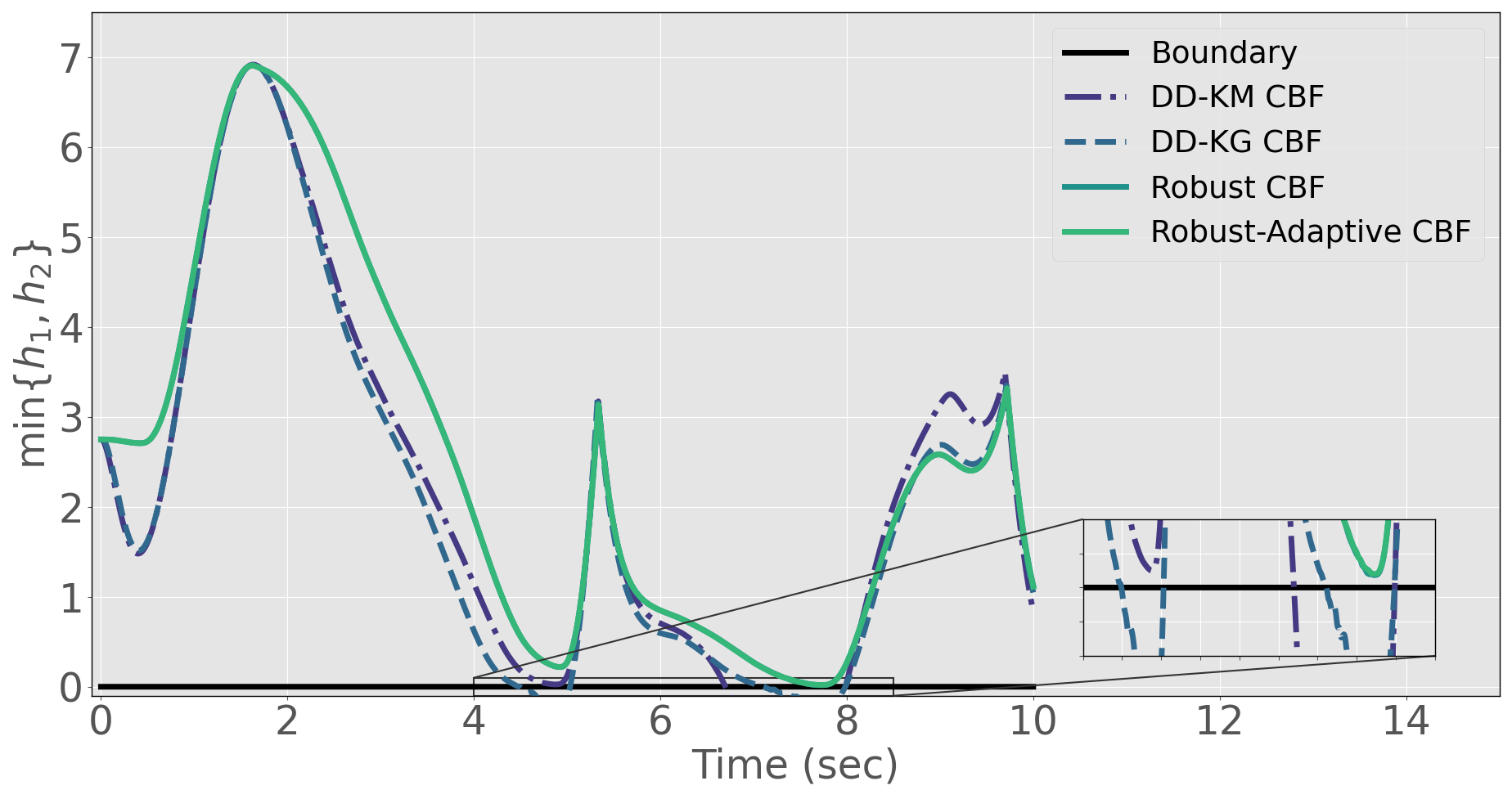}
        \caption{\small{Evolutions of $h_1$ and $h_2$ for the various controllers considered in the double-integrator example.}}
        \label{fig.cbf_single_integrator}
    \end{center}
\end{figure}

\begin{figure}
    \begin{center}
        \includegraphics[width=1\columnwidth,clip]{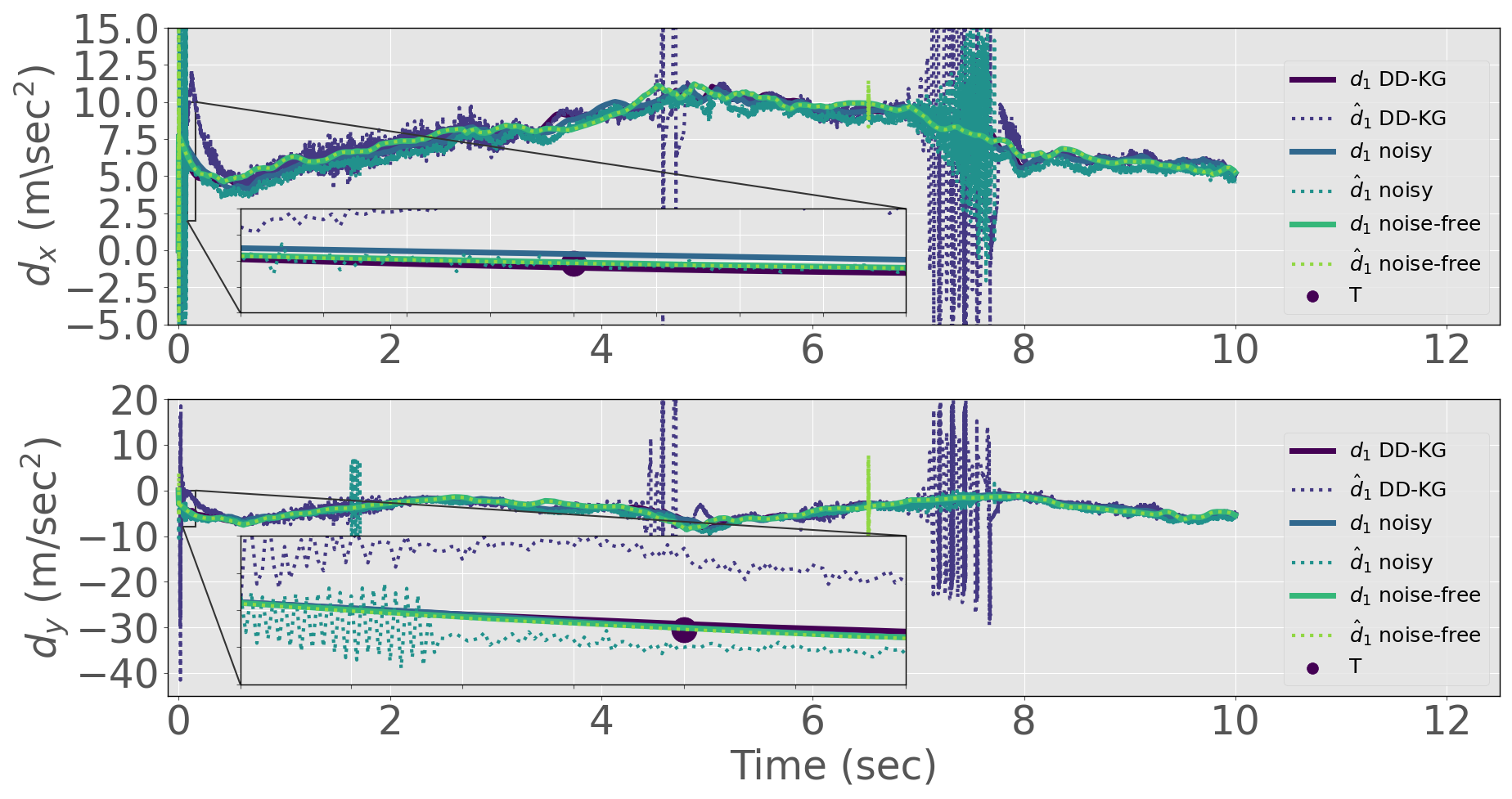}
        \caption{\small{The estimates $\hat d_x$, $\hat d_y$ of the unknown wind gusts ($d_x$ and $d_x$). In our scheme, the estimates converge to the true values within the fixed-time $T$ without noise, and converge to a close approximation in the presence of measurement noise.}}
        \label{fig.estimation_single_integrator}
    \end{center}
\end{figure}

\begin{figure}
    \begin{center}
        \includegraphics[width=1\columnwidth,clip]{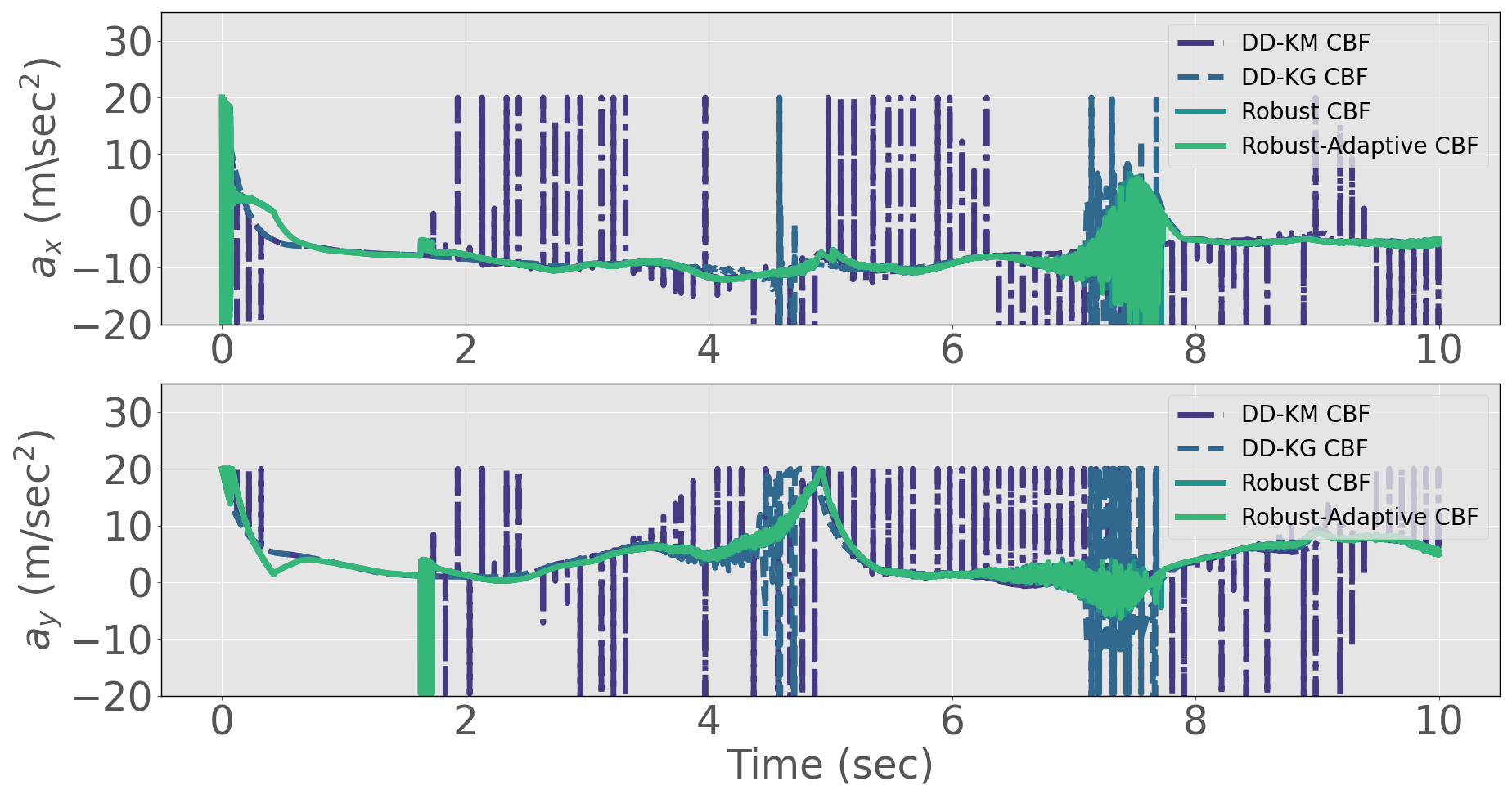}
        \caption{\small{Control inputs for the double-integrator example.}}
        \label{fig.u_single_integrator}
    \end{center}
\end{figure}

The resulting paths taken by the simulated CBF-controlled vehicles (Koopman-based naive, robust, and robust-adaptive), as well as the path taken for the nominally controlled vehicle without disturbance estimation are displayed in Figure \ref{fig.xy_single_integrator}. Here, only the robust and robust-adaptive CBF controllers that use our fixed-time identification approach preserve safety (as seen in Figure \ref{fig.cbf_single_integrator}). As the data-driven Koopman matrix (\cite{Bruder2021DataDriven}) and generator (\cite{Klus2020KoopmanGenerator}) approaches are non-recursive and unable to quantify the identification error, they are neither sufficiently responsive nor accurate enough to guarantee safety in this example. Figure \ref{fig.estimation_single_integrator} highlights that our disturbance estimates indeed converge to the true values within the fixed-time $T=0.12$ sec, computed using \eqref{eq.koopman_generator_settling_time}, and the control inputs are shown in Figure \ref{fig.u_single_integrator}. We further note that even when measurement noise is injected into the system, the adaptation-based approach succeeds in both reconstructing the unknown disturbance to within a small error and preserving safety. We leave quantification of this measurement error and any error associated with representing the infinite-dimensional Koopman operator in a finite-dimensional subspace to future work.



\section{Conclusion}\label{sec.conclusion}
We introduced a safe control synthesis using Koopman-based fixed-time system identification. We showed that under mild assumptions we can learn the unknown, additive, nonlinear vector field perturbing the system dynamics within a fixed-time independent of the initial estimate. The a priori knowledge of this identification guarantee allows us to derive robust and robust-adaptive control barrier function conditions suitable for use in a standard class of quadratic program based controllers.

We recognize that there are practical limitations to our method, including the need to measure the state derivative and to be able to exactly represent the linear, infinite-dimensional Koopman dynamical system with a finite-rank operator. Though we demonstrated some robustness to measurement noise in our simulated study, in the future we will seek to relax these assumptions by analyzing the use of observers and filters for state and state-derivative approximation and by seeking to quantify the residual error associated with projecting the infinite-dimensional Koopman operator onto a finite-dimensional subspace.






\begin{ack}
The authors would like to acknowledge the support of the National Science Foundation (NSF) through grants 1931982 and 1942907.
\end{ack}

\bibliography{ifacconf}             

\begin{thebibliography}{4}
\providecommand{\natexlab}[1]{#1}
\providecommand{\url}[1]{\texttt{#1}}
\providecommand{\urlprefix}{URL }
\expandafter\ifx\csname urlstyle\endcsname\relax
  \providecommand{\doi}[1]{doi:\discretionary{}{}{}#1}\else
  \providecommand{\doi}{doi:\discretionary{}{}{}\begingroup
  \urlstyle{rm}\Url}\fi

\bibitem[{Able(1956)}]{Abl:56}
Able, B. (1956).
\newblock Nucleic acid content of microscope.
\newblock \emph{Nature}, 135, 7--9.

\bibitem[{Able et~al.(1954)Able, Tagg, and Rush}]{AbTaRu:54}
Able, B., Tagg, R., and Rush, M. (1954).
\newblock Enzyme-catalyzed cellular transanimations.
\newblock In A.~Round (ed.), \emph{Advances in Enzymology}, volume~2, 125--247.
  Academic Press, New York, 3rd edition.

\bibitem[{Keohane(1958)}]{Keo:58}
Keohane, R. (1958).
\newblock \emph{Power and Interdependence: World Politics in Transitions}.
\newblock Little, Brown \& Co., Boston.

\bibitem[{Powers(1985)}]{Pow:85}
Powers, T. (1985).
\newblock Is there a way out?
\newblock \emph{Harpers}, 35--47.

\end{thebibliography}


\begin{thebibliography}{28}
\providecommand{\natexlab}[1]{#1}
\providecommand{\url}[1]{\texttt{#1}}
\providecommand{\urlprefix}{URL }
\expandafter\ifx\csname urlstyle\endcsname\relax
  \providecommand{\doi}[1]{doi:\discretionary{}{}{}#1}\else
  \providecommand{\doi}{doi:\discretionary{}{}{}\begingroup
  \urlstyle{rm}\Url}\fi

\bibitem[{Ames et~al.(2017)Ames, Xu, Grizzle, and Tabuada}]{ames2017control}
Ames, A.D., Xu, X., Grizzle, J.W., and Tabuada, P. (2017).
\newblock Control barrier function based quadratic programs for safety critical
  systems.
\newblock \emph{IEEE Trans. on Automatic Control}, 62(8), 3861--3876.

\bibitem[{Black et~al.(2022{\natexlab{a}})Black, Arabi, and
  Panagou}]{Black2022Fixed}
Black, M., Arabi, E., and Panagou, D. (2022{\natexlab{a}}).
\newblock Fixed-time parameter adaptation for safe control synthesis.
\newblock \emph{arXiv preprint arXiv:2204.10453}.

\bibitem[{Black et~al.(2020)Black, Garg, and Panagou}]{Black2020Quadratic}
Black, M., Garg, K., and Panagou, D. (2020).
\newblock A quadratic program based control synthesis under spatiotemporal
  constraints and non-vanishing disturbances.
\newblock In \emph{2020 59th IEEE Conference on Decision and Control (CDC)},
  2726--2731.

\bibitem[{Black et~al.(2022{\natexlab{b}})Black, Jankovic, Sharma, and
  Panagou}]{black2022intersection}
Black, M., Jankovic, M., Sharma, A., and Panagou, D. (2022{\natexlab{b}}).
\newblock Future-focused control barrier functions for autonomous vehicle
  control.
\newblock \emph{arXiv preprint arXiv:2204.00127}.

\bibitem[{Blanchini(1999)}]{blanchini1999set}
Blanchini, F. (1999).
\newblock Set invariance in control.
\newblock \emph{Automatica}, 35(11), 1747--1767.

\bibitem[{Bruder et~al.(2021)Bruder, Fu, Gillespie, Remy, and
  Vasudevan}]{Bruder2021DataDriven}
Bruder, D., Fu, X., Gillespie, R.B., Remy, C.D., and Vasudevan, R. (2021).
\newblock Data-driven control of soft robots using koopman operator theory.
\newblock \emph{IEEE Transactions on Robotics}, 37(3), 948--961.
\newblock \doi{10.1109/TRO.2020.3038693}.

\bibitem[{Brunton et~al.(2016)Brunton, Brunton, Proctor, and
  Kutz}]{brunton2016koopman}
Brunton, S.L., Brunton, B.W., Proctor, J.L., and Kutz, J.N. (2016).
\newblock Koopman invariant subspaces and finite linear representations of
  nonlinear dynamical systems for control.
\newblock \emph{PloS one}, 11(2), e0150171.

\bibitem[{Davoudi et~al.(2020)Davoudi, Taheri, Duraisamy, Jayaraman, and
  Kolmanovsky}]{davoudi2020quad}
Davoudi, B., Taheri, E., Duraisamy, K., Jayaraman, B., and Kolmanovsky, I.
  (2020).
\newblock Quad-rotor flight simulation in realistic atmospheric conditions.
\newblock \emph{AIAA Journal}, 58(5), 1992--2004.

\bibitem[{Drmač et~al.(2021)Drmač, Mezić, and
  Mohr}]{Drmac2021KoopmanGenerator}
Drmač, Z., Mezić, I., and Mohr, R. (2021).
\newblock Identification of nonlinear systems using the infinitesimal generator
  of the koopman semigroup—a numerical implementation of the
  mauroy–goncalves method.
\newblock \emph{Mathematics}, 9(17), 2075.
\newblock \doi{10.3390/math9172075}.

\bibitem[{Folkestad et~al.(2020)Folkestad, Chen, Ames, and
  Burdick}]{folkestad2020data}
Folkestad, C., Chen, Y., Ames, A.D., and Burdick, J.W. (2020).
\newblock Data-driven safety-critical control: Synthesizing control barrier
  functions with koopman operators.
\newblock \emph{IEEE Control Systems Letters}, 5(6), 2012--2017.

\bibitem[{Frigola and Rasmussen(2013)}]{frigola2013integrated}
Frigola, R. and Rasmussen, C.E. (2013).
\newblock Integrated pre-processing for bayesian nonlinear system
  identification with gaussian processes.
\newblock In \emph{52nd IEEE Conference on Decision and Control}, 5371--5376.
  IEEE.

\bibitem[{Gonzalez and Yu(2018)}]{gonzalez2018non}
Gonzalez, J. and Yu, W. (2018).
\newblock Non-linear system modeling using lstm neural networks.
\newblock \emph{IFAC-PapersOnLine}, 51(13), 485--489.

\bibitem[{Haseli and Cortés(2021)}]{Haseli2021Approximation}
Haseli, M. and Cortés, J. (2021).
\newblock Data-driven approximation of koopman-invariant subspaces with tunable
  accuracy.
\newblock In \emph{2021 American Control Conference (ACC)}, 470--475.
\newblock \doi{10.23919/ACC50511.2021.9483259}.

\bibitem[{Jankovic(2018)}]{Jankovic2018Robust}
Jankovic, M. (2018).
\newblock Robust control barrier functions for constrained stabilization of
  nonlinear systems.
\newblock \emph{Automatica}, 96, 359--367.

\bibitem[{Klus et~al.(2020)Klus, Nüske, Peitz, Niemann, Clementi, and
  Schütte}]{Klus2020KoopmanGenerator}
Klus, S., Nüske, F., Peitz, S., Niemann, J.H., Clementi, C., and Schütte, C.
  (2020).
\newblock Data-driven approximation of the koopman generator: Model reduction,
  system identification, and control.
\newblock \emph{Physica D: Nonlinear Phenomena}, 406, 132416.
\newblock \doi{https://doi.org/10.1016/j.physd.2020.132416}.

\bibitem[{Lopez et~al.(2021)Lopez, Slotine, and How}]{Lopez2020Robust}
Lopez, B.T., Slotine, J.J.E., and How, J.P. (2021).
\newblock Robust adaptive control barrier functions: An adaptive and
  data-driven approach to safety.
\newblock \emph{IEEE Control Systems Letters}, 5(3), 1031--1036.
\newblock \doi{10.1109/LCSYS.2020.3005923}.

\bibitem[{Mauroy and Goncalves(2020)}]{Mauroy2020Koopman}
Mauroy, A. and Goncalves, J. (2020).
\newblock Koopman-based lifting techniques for nonlinear systems
  identification.
\newblock \emph{IEEE Transactions on Automatic Control}, 65(6), 2550--2565.
\newblock \doi{10.1109/TAC.2019.2941433}.

\bibitem[{Mauroy et~al.(2020)Mauroy, Susuki, and
  Mezi{\'c}}]{mauroy2020koopmanbook}
Mauroy, A., Susuki, Y., and Mezi{\'c}, I. (2020).
\newblock \emph{Koopman operator in systems and control}.
\newblock Springer.

\bibitem[{Ortega et~al.(2022)Ortega, Bobtsov, and
  Nikolaev}]{Ortega2022Parameter}
Ortega, R., Bobtsov, A., and Nikolaev, N. (2022).
\newblock Parameter identification with finite-convergence time alertness
  preservation.
\newblock \emph{IEEE Control Systems Letters}, 6, 205--210.
\newblock \doi{10.1109/LCSYS.2021.3057012}.

\bibitem[{Polyakov(2012)}]{polyakov2012nonlinear}
Polyakov, A. (2012).
\newblock Nonlinear feedback design for fixed-time stabilization of linear
  control systems.
\newblock \emph{IEEE Transactions on Automatic Control}, 57(8), 2106.

\bibitem[{R{\'\i}os et~al.(2017)R{\'\i}os, Efimov, Moreno, Perruquetti, and
  Rueda-Escobedo}]{rios2017time}
R{\'\i}os, H., Efimov, D., Moreno, J.A., Perruquetti, W., and Rueda-Escobedo,
  J.G. (2017).
\newblock Time-varying parameter identification algorithms: Finite and
  fixed-time convergence.
\newblock \emph{IEEE Transactions on Automatic Control}, 62(7), 3671--3678.

\bibitem[{Schoellig et~al.(2012)Schoellig, Wiltsche, and
  D'Andrea}]{Schoellig2012PeriodicQuadrotor}
Schoellig, A.P., Wiltsche, C., and D'Andrea, R. (2012).
\newblock Feed-forward parameter identification for precise periodic
  quadrocopter motions.
\newblock In \emph{2012 American Control Conference (ACC)}, 4313--4318.
\newblock \doi{10.1109/ACC.2012.6315248}.

\bibitem[{Taylor and Ames(2020)}]{Taylor2020Adaptive}
Taylor, A.J. and Ames, A.D. (2020).
\newblock Adaptive safety with control barrier functions.
\newblock In \emph{2020 American Control Conference (ACC)}, 1399--1405.
\newblock \doi{10.23919/ACC45564.2020.9147463}.

\bibitem[{Wang et~al.(2022)Wang, Lyu, Wen, Shi, Zhu, and
  Huang}]{Wang2022Robust}
Wang, S., Lyu, B., Wen, S., Shi, K., Zhu, S., and Huang, T. (2022).
\newblock Robust adaptive safety-critical control for unknown systems with
  finite-time elementwise parameter estimation.
\newblock \emph{IEEE Transactions on Systems, Man, and Cybernetics: Systems},
  1--11.
\newblock \doi{10.1109/TSMC.2022.3203176}.

\bibitem[{Williams et~al.(2015)Williams, Kevrekidis, and
  Rowley}]{williams2015data}
Williams, M.O., Kevrekidis, I.G., and Rowley, C.W. (2015).
\newblock A data--driven approximation of the koopman operator: Extending
  dynamic mode decomposition.
\newblock \emph{Journal of Nonlinear Science}, 25(6), 1307--1346.

\bibitem[{Zancato and Chiuso(2021)}]{zancato2021novel}
Zancato, L. and Chiuso, A. (2021).
\newblock A novel deep neural network architecture for non-linear system
  identification.
\newblock \emph{IFAC-PapersOnLine}, 54(7), 186--191.

\bibitem[{Zhou and Schwager(2014)}]{Zhou2014vector}
Zhou, D. and Schwager, M. (2014).
\newblock Vector field following for quadrotors using differential flatness.
\newblock In \emph{2014 IEEE International Conference on Robotics and
  Automation (ICRA)}, 6567--6572.
\newblock \doi{10.1109/ICRA.2014.6907828}.

\bibitem[{Zinage and Bakolas(2022)}]{zinage2022neural}
Zinage, V. and Bakolas, E. (2022).
\newblock Neural koopman control barrier functions for safety-critical control
  of unknown nonlinear systems.
\newblock \emph{arXiv preprint arXiv:2209.07685}.

\end{thebibliography}

\appendix

\section{Proof of Corollary \ref{cor.disturbance_fixed_time}}\label{app.proof_d_error}
Using \eqref{eq.dhat} and \eqref{eq.w_mat} we can express the disturbance vector field error as
\begin{align}
    \Tilde{d}(\bb{x}(t)) &= d(\bb{x}(t)) - \hat d(\bb{x}(t)), \nonumber \\ &= \bb{W}(t)\Tilde{\bb{\lambda}}(t), \nonumber
\end{align}
and thus can observe that $\|d(\bb{x}(t))\|_{\infty} = \|\bb{W}(t)\Tilde{\bb{\lambda}}(t)\|_{\infty} \leq \sigma_{max}(\bb{W}(t))\|\Tilde{\bb{\lambda}}(t)\|_{\infty}$. Then, using \cite[Corollary 1]{Black2022Fixed} we obtain that
$\|\Tilde{\bb{\lambda}}(t)\|_{\infty} \leq \Lambda \tan^{\frac{w}{2}}(A(t))$ for all $t \leq T$, where $\Lambda$, $A(t)$, and $T$ are given by \eqref{eq.Lambda_gain}, \eqref{eq.A(t)}, \eqref{eq.koopman_generator_settling_time} respectively, and $\|\Tilde{\bb{\lambda}}(t)\|_{\infty} = 0$ for all $t > T$.

Then, to obtain the $\Xi$ term in \eqref{eq.xi_term}, observe that with $\hat{\bb{\lambda}}(0) = \mathbf{0}_{N^2 \times 1}$ and the assumption that $\|d(\bb{x})\|_{\infty} \leq D$, $\forall \bb{x} \in \mathcal{X}$, it follows that at $t = 0$
\begin{equation}\nonumber
    \sigma_{min}(\bb{W})\|\Tilde{\bb{\lambda}}\|_{\infty} \leq \|\bb{W}\Tilde{\bb{\lambda}}\|_{\infty} = \|\Tilde{d}(\bb{x})\|_{\infty} \leq 2D,
\end{equation}
from which we obtain that
\begin{equation}\nonumber
    \|\Tilde{\bb{\lambda}}(0)\|_{\infty} \leq \frac{2D}{\sigma_{min}\left(\bb{W}(0)\right)}.
\end{equation}
Thus we obtain $\bb{l} = \frac{2D}{\sigma_{min}\left(\bb{W}(0)\right)} \cdot \mathbf{1}_{N^2 \times 1}$, and this completes the proof.

\end{document}